# On a certain relation between Legendre's conjecture and Bertrand's postulate

Tsutomu Hashimoto

July, 2008

## Abstract

Given a natural number $n$, if $\pi$ is the prime-counting function, then

$$\pi((n+1)^2) - \pi(n^2) = \pi(2n) - \pi(n) + 1 - \phi_T(n^2, 2n, \pi(n)),$$

where $\phi_T$ is a certain transformation of Legendre's function.

## Introduction

Let $\pi(x)$ be the prime-counting function [3].

Legendre's conjecture [2] is

"There's always a prime between $n^2$ and $(n+1)^2$."

That is, $\pi((n+1)^2) - \pi(n^2) \geq 1$.

Bertrand's postulate [1] is

"There's always a prime between $n$ and $2n$."

That is, $\pi(2n) - \pi(n) \geq 1$.

In this paper, we prove a certain relation between Legendre's conjecture and Bertrand's postulate in terms of the following transformation of Legendre's function $\phi$ [4].

Given natural numbers $M_1$, $M_2$, and $n$, define the function

$$\phi_T(M_1, M_2, \pi(n)) := \sum_{p_i \leq n} \left\lfloor \frac{M_1 \bmod p_i + M_2 \bmod p_i}{p_i} \right\rfloor - \sum_{p_i < p_j \leq n} \left\lfloor \frac{M_1 \bmod p_i p_j + M_2 \bmod p_i p_j}{p_i p_j} \right\rfloor$$
$$+ \sum_{p_i < p_j < p_k \leq n} \left\lfloor \frac{M_1 \bmod p_i p_j p_k + M_2 \bmod p_i p_j p_k}{p_i p_j p_k} \right\rfloor - \cdots,$$

where the numbers $p_i, p_j, p_k, \cdots$ are the primes less than or equal to $n$, and $\lfloor x \rfloor$ is the floor function.

**Theorem 1.**

*Given a natural number $n$, if $\pi$ is the prime-counting function, then*
$$\pi((n+1)^2) - \pi(n^2) = \pi(2n) - \pi(n) + 1 - \phi_T(n^2, 2n, \pi(n)).$$

For example, take $n = 6$:

A) $\pi(7^2) - \pi(6^2) = \pi(2 \cdot 6) - \pi(6) + 1 - \phi_T(6^2, 2 \cdot 6, \pi(6))$

B) $\pi(49) - \pi(36) = \pi(12) - \pi(6) + 1$

$$-\left(\left\lfloor\frac{36 \bmod 2 \;+\; 12 \bmod 2}{2}\right\rfloor + \left\lfloor\frac{36 \bmod 3 \;+\; 12 \bmod 3}{3}\right\rfloor + \left\lfloor\frac{36 \bmod 5 \;+\; 12 \bmod 5}{5}\right\rfloor - \left\lfloor\frac{36 \bmod 6 \;+\; 12 \bmod 6}{6}\right\rfloor\right.$$

$$\left. - \left\lfloor\frac{36 \bmod 10 \;+\; 12 \bmod 10}{10}\right\rfloor - \left\lfloor\frac{36 \bmod 15 \;+\; 12 \bmod 15}{15}\right\rfloor + \left\lfloor\frac{36 \bmod 30 \;+\; 12 \bmod 30}{30}\right\rfloor\right)$$

C) $15 - 11 = 5 - 3 + 1 - (0 + 0 + 0 - 0 - 0 - 1 + 0)$

D) $4 = 4$

**Corollary 1.**

*Legendre's conjecture is true if and only if*
$$\phi_T(n^2, 2n, \pi(n)) \leq \pi(2n) - \pi(n) \text{ for all } n.$$
*In particular, for fixed $n$, if $\phi_T(n^2, 2n, \pi(n)) \leq 1$, then $\pi((n+1)^2) - \pi(n^2) \geq 1$.*
*More generally, if $\phi_T(n^2, 2n, \pi(n))$ is bounded above, then*
$$\pi((n+1)^2) - \pi(n^2) \to \infty \text{ as } n \to \infty.$$

*Proof.* The first two statements follow from Theorem 1 and Bertrand's postulate. The third statement follows from Theorem 1, using the fact that $\pi(2n) - \pi(n) \to \infty$ as $n \to \infty$, which Ramanujan proved in [6] (see also [1]). ∎

Define the functions $B_t(n) := \pi(2n) - \pi(n)$, and $P_t(n) := \phi_T(n^2, 2n, \pi(n))$.

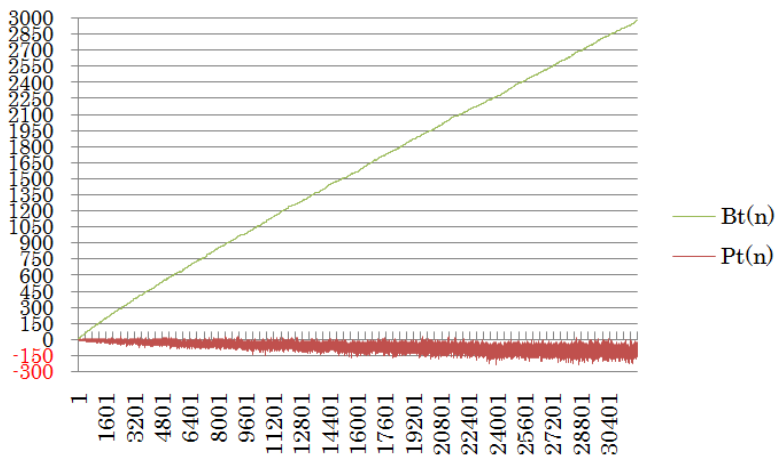

# 1. $T_{\sqrt{N}}(M)$ (Multiple-counting function)

Given a positive integer $N$, let $M$ denote a natural number such that $\lfloor\sqrt{N}\rfloor \leq M < (\lfloor\sqrt{N}\rfloor + 1)^2$. Define the multiple-counting function

$$T_{\sqrt{N}}(M) := M - \phi\left(M, \pi(\sqrt{N})\right),$$

where $\phi$ is Legendre's function [4].
Then

$$T_{\sqrt{N}}(M) = \sum_{p_i \leq \sqrt{N}} \left\lfloor \frac{M}{p_i} \right\rfloor - \sum_{p_i < p_j \leq \sqrt{N}} \left\lfloor \frac{M}{p_i p_j} \right\rfloor + \sum_{p_i < p_j < p_k \leq \sqrt{N}} \left\lfloor \frac{M}{p_i p_j p_k} \right\rfloor - \cdots,$$

where the numbers $p_i, p_j, p_k, \cdots$ are the primes less than or equal to the square root of $N$ [5].

Let $B_{\sqrt{N}}$ be the set of all "denominators" of $T_{\sqrt{N}}(M)$ (e.g., $p_1$, $p_1 p_2$, $p_1 p_2 p_3$, etc.) and define

$$\sum_{\beta_{\sqrt{N}}}^{+-} \eta(\beta_{\sqrt{N}}) := \sum_{p_i \leq \sqrt{N}} \eta(p_i) - \sum_{p_i < p_j \leq \sqrt{N}} \eta(p_i p_j) + \sum_{p_i < p_j < p_k \leq \sqrt{N}} \eta(p_i p_j p_k) - \cdots,$$

where $\beta_{\sqrt{N}}$ runs over all elements of $B_{\sqrt{N}}$ and $\eta$ is an arbitrary function of $\beta_{\sqrt{N}}$.
Then

$$\sum_{\beta_{\sqrt{N}}}^{+-} \eta(\beta_{\sqrt{N}}) = T_{\sqrt{N}}(M),$$

where $\eta(\beta_{\sqrt{N}}) = \left\lfloor \frac{M}{\beta_{\sqrt{N}}} \right\rfloor$.

2. $T_n(n^2 + 2n) - T_n(n^2) = T_n(2n) + \sum_{\beta_n}^{+-} \left[ \frac{H_{\beta_n}^{n^2} + H_{\beta_n}^{2n}}{\beta_n} \right]$

Given a natural number $n$, if

$$\begin{cases} H_{\beta_n}^{n^2+2n} := (n^2 + 2n) \bmod \beta_n \; ; \; 0 \leq H_{\beta_n}^{n^2+2n} < \beta_n, \\ n^2 + 2n = \left\lfloor \dfrac{n^2 + 2n}{\beta_n} \right\rfloor \cdot \beta_n + H_{\beta_n}^{n^2+2n} \end{cases}$$

and if

$$\begin{cases} H_{\beta_n}^{n^2} := n^2 \bmod \beta_n \; ; \; 0 \leq H_{\beta_n}^{n^2} < \beta_n, \\ n^2 = \left\lfloor \dfrac{n^2}{\beta_n} \right\rfloor \cdot \beta_n + H_{\beta_n}^{n^2}, \end{cases}$$

then

$$T_n(n^2 + 2n) - T_n(n^2) = \sum_{\beta_n}^{+-} \left( \left\lfloor \frac{n^2 + 2n}{\beta_n} \right\rfloor - \left\lfloor \frac{n^2}{\beta_n} \right\rfloor \right) = \sum_{\beta_n}^{+-} \left( \frac{n^2 + 2n - H_{\beta_n}^{n^2+2n}}{\beta_n} - \frac{n^2 - H_{\beta_n}^{n^2}}{\beta_n} \right)$$

$$= \sum_{\beta_n}^{+-} \left( \frac{2n - H_{\beta_n}^{n^2+2n} + H_{\beta_n}^{n^2}}{\beta_n} \right).$$

On the other hand, if

$$\begin{cases} H_{\beta_n}^{2n} := 2n \bmod \beta_n \; ; \; 0 \leq H_{\beta_n}^{2n} < \beta_n, \\ 2n = \left\lfloor \dfrac{2n}{\beta_n} \right\rfloor \cdot \beta_n + H_{\beta_n}^{2n}, \end{cases}$$

then

$$T_n(2n) + \sum_{\beta_n}^{+-} \left[ \frac{H_{\beta_n}^{n^2} + H_{\beta_n}^{2n}}{\beta_n} \right] = \sum_{\beta_n}^{+-} \left( \left\lfloor \frac{2n}{\beta_n} \right\rfloor + \left\lfloor \frac{H_{\beta_n}^{n^2} + H_{\beta_n}^{2n}}{\beta_n} \right\rfloor \right)$$

$$= \sum_{\beta_n}^{+-} \left( \frac{2n - H_{\beta_n}^{2n}}{\beta_n} + \frac{H_{\beta_n}^{n^2} + H_{\beta_n}^{2n} - H_{\beta_n}^{n^2+2n}}{\beta_n} \right) = \sum_{\beta_n}^{+-} \left( \frac{2n - H_{\beta_n}^{n^2+2n} + H_{\beta_n}^{n^2}}{\beta_n} \right).$$

Therefore

$$T_n(n^2 + 2n) - T_n(n^2) = T_n(2n) + \sum_{\beta_n}^{+-} \left[ \frac{H_{\beta_n}^{n^2} + H_{\beta_n}^{2n}}{\beta_n} \right].$$

## 3. Proof of Theorem 1

Recall the definition of the function

$$\phi_T\left(M_1, M_2, \pi(\sqrt{N})\right) := \sum_{\beta_{\sqrt{N}}}^{+-} \eta'(\beta_{\sqrt{N}}),$$

where $\eta'(\beta_{\sqrt{N}}) = \left\lfloor \dfrac{M_1 \bmod \beta_{\sqrt{N}} + M_2 \bmod \beta_{\sqrt{N}}}{\beta_{\sqrt{N}}} \right\rfloor$. (cf. section 1)

Legendre's formula [4] is $\phi\left(M, \pi(\sqrt{N})\right) = \pi(M) - \pi(\sqrt{N}) + 1$.

From section 2, we have the following equation:

$$T_n(n^2 + 2n) - T_n(n^2) = T_n(2n) + \sum_{\beta_n}^{+-} \left\lfloor \frac{H_{\beta_n}^{n^2} + H_{\beta_n}^{2n}}{\beta_n} \right\rfloor.$$

This is equivalent to the equation

$$2n - \left(\pi(n^2 + 2n) - \pi(n^2)\right) = 2n - (\pi(2n) - \pi(n) + 1) + \phi_T\left(n^2, 2n, \pi(n)\right),$$

because $T_n(M) = M - \phi(M, \pi(n)) = M - (\pi(M) - \pi(n) + 1)$ and

$$\left\lfloor \frac{H_{\beta_n}^{n^2} + H_{\beta_n}^{2n}}{\beta_n} \right\rfloor = \left\lfloor \frac{n^2 \bmod \beta_n + 2n \bmod \beta_n}{\beta_n} \right\rfloor.$$

It follows, since $\pi((n+1)^2) = \pi(n^2 + 2n)$, that

$$\pi((n+1)^2) - \pi(n^2) = \pi(2n) - \pi(n) + 1 - \phi_T\left(n^2, 2n, \pi(n)\right).$$

∎

### Acknowledgement

I thank Jonathan Sondow for several suggestions which improved the paper.

Shiga 520-2412 JAPAN

t-hashimoto@aquablue.ne.jp